\documentclass[aop,MSNbibl,dvips]{arximspdf}
\usepackage{mathbh}
\usepackage{graphicx}

\doi{10.1214/13-AOP847} 
\volume{42}
\issue{4}
\pubyear{2014}
\firstpage{1374}
\lastpage{1395}

\makeatletter
\newtheorem{theorem}{Theorem}[section]
\newtheorem{lemma}[theorem]{Lemma}
\newtheorem{corollary}[theorem]{Corollary}
\newtheorem{proposition}[theorem]{Proposition}
\newproclaim{definition}[theorem]{Definition}
\newproclaim{remark}[theorem]{Remark}
\newproclaim{example}[theorem]{Example}
\makeatother

\begin{document}
\begin{frontmatter}

\title{Asymptotic distribution of complex zeros of~random analytic functions}
\runtitle{Complex zeros of random analytic functions}

\begin{aug}
\author[A]{\fnms{Zakhar} \snm{Kabluchko}\corref{}\ead[label=e1]{zakhar.kabluchko@uni-ulm.de}}
\and
\author[B]{\fnms{Dmitry} \snm{Zaporozhets}\ead[label=e2]{zap1979@gmail.com}}
\runauthor{Z. Kabluchko and D. Zaporozhets}
\affiliation{Ulm University and Steklov Institute of Mathematics}
\address[A]{Institute of Stochastics\\
Ulm University\\
Helmholtzstr. 18\\
89069 Ulm\\
Germany\\
\printead{e1}}
\address[B]{St. Petersburg Branch\\ 
Steklov Institute of Mathematics\\
Fontanka Str. 27\\
191011 St. Petersburg\\
Russia\\
\printead{e2}}
\end{aug}

\received{\smonth{5} \syear{2012}}
\revised{\smonth{1} \syear{2013}}

%
\begin{abstract}
Let $\xi_0,\xi_1,\ldots$ be independent identically distributed
complex-\break valued random variables such that $\mathbb{E}\log(1+|\xi
_0|)<\infty$.
We consider random analytic functions of the form
\[
\mathbf{G}_n(z)=\sum_{k=0}^{\infty}
\xi_k f_{k,n} z^k,
\]
where $f_{k,n}$ are deterministic complex coefficients. Let $\mu_n$ be
the random measure counting the complex zeros of $\mathbf{G}_n$
according to
their multiplicities. Assuming essentially that $-\frac1n \log
f_{[tn], n}\to u(t)$ as $n\to\infty$, where $u(t)$ is some function,
we show that the measure $\frac1n \mu_n$ converges in probability to
some deterministic measure $\mu$ which is characterized in terms of
the Legendre--Fenchel transform of $u$. The limiting measure $\mu$
does not depend on the distribution of the $\xi_k$'s. This result is
applied to several ensembles of random analytic functions including the
ensembles corresponding to the three two-dimensional geometries of
constant curvature. As another application, we prove a random
polynomial analogue of the circular law for random matrices.
\end{abstract}

%

\end{frontmatter}

\section{Introduction}\label{secintro}
\subsection{Statement of the problem}\label{subsecintro}
Let $\xi_0,\xi_1,\ldots$ be nondegenerate independent identically
distributed (i.i.d.) random variables with complex values. The simplest
ensemble of random polynomials are the \textit{Kac polynomials} defined
as
\[
\mathbf{K}_n(z)=\sum_{k=0}^n
\xi_k z^k.
\]
The distribution of zeros of Kac polynomials has been much studied; see
\cite{hammersley,sparosur,arnold,sheppvanderbei,ibrzeit,shmerlinghochberg,hughesnikeghbali,izlog}.
It is known that under a very mild moment assumption, the complex zeros
of $\mathbf{K}_n$ cluster asymptotically near the unit circle
$\mathbb{T}=\{|z|=1\}$ and that the distribution of zeros is
asymptotically uniform with regard to the argument. To make this
precise, we need to introduce some notation. Let $G$ be an analytic
function in some domain $D\subset\mathbb{C}$. Assuming that $G$ does
not vanish identically,\vadjust{\goodbreak} we consider a measure
$\mu_{G}$ counting the complex zeros of $G$ according to their
multiplicities:\vspace*{-3pt}
\[
\mu_{G}=\sum_{z\in D\dvtx  G(z)=0} n_G(z)
\delta(z).
\]
%
Here, $n_{G}(z)$ is the multiplicity of the zero at $z$ and $\delta(z)$
is the unit point mass at $z$. If $G$ vanishes identically, we put
$\mu_G=0$. Then, Ibragimov and Zaporozhets \cite{izlog} proved that the
following two conditions are equivalent:\vspace*{-2.5pt}
\begin{longlist}[(2)]
\item[(1)] With probability $1$, the sequence of measures $\frac1n
\mu_{\mathbf{K}_n}$ converges as $n\to\infty$ weakly to the uniform
probability distribution on $\mathbb{T}$.

\item[(2)] $\mathbb{E}\log(1+|\xi_0|)<\infty$.\vspace*{-1.5pt}
\end{longlist}

\begin{figure}

\includegraphics{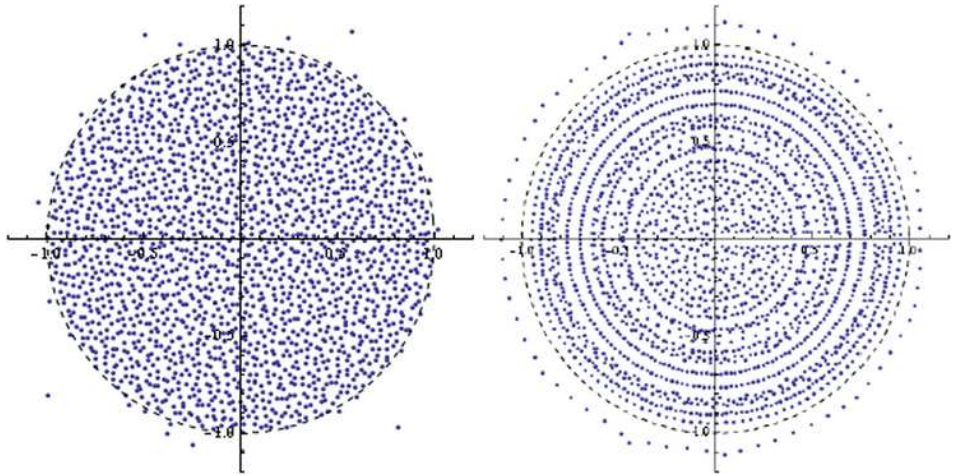}

\caption{Zeros\vspace*{-1pt} of the Weyl random polynomial $\mathbf{W}_n(z)=\sum_{k=0}^n
\xi_k \frac{z^k}{\sqrt{k!}}$ of degree $n=2000$. The zeros were divided
by $\sqrt n$. Left: Complex normal coefficients. Right: Coefficients
are positive with $\mathbb{P}[\log\xi_k>t]=t^{-4}$ for $t>1$.
In both cases, the limiting distribution of zeros is uniform on the
unit~disk.}\label{figweyluniversality}
\end{figure}

Along with the Kac polynomials, many other remarkable ensembles of
random polynomials (or, more generally, random power series) appeared
in the literature. These ensembles are usually characterized by
invariance properties with respect to certain groups of transformations
and have the general form\vspace*{-1pt}
\[
\mathbf{G}_n(z)=\sum_{k=0}^{\infty}
\xi_k f_{k,n}z^k,
\]
where $\xi_0,\xi_1,\ldots$ are i.i.d. complex-valued random variables
and $f_{k,n}$ are complex deterministic coefficients.
The aim of the present work is to study the distribution of zeros of
$\mathbf{G}_n$ asymptotically as $n\to\infty$. We will show that under
certain assumptions on the coefficients $f_{k,n}$, the random measure
$\frac1n \mu_{\mathbf{G}_n}$ converges, as $n\to\infty$, to~some
limiting deterministic measure $\mu$. The limiting measure $\mu$ does
not depend on the distribution of the random variables $\xi_k$; see
Figure~\ref{figweyluniversality}.
Results of this\vadjust{\goodbreak} type are known in the context of
random matrices; see, for example,~\cite{taovu}.
However, the literature on random polynomials and random analytic
functions usually concentrates on the Gaussian case, since in this case
explicit calculations are possible; see, for example,
\cite{hammersley,edelmankostlan,peresbook,sodintsirelson,sheppvanderbei,schiffmanzelditchbundles,shiffmanzelditch1,bloomshiffman,sodinECM,farahmandbook,bharuchareidbook}.
The only ensemble of random polynomials for which the independence of
the limiting distribution of zeros on the distribution of the
coefficients is well understood is the Kac ensemble;
see~\cite{sparosur,arnold,ibrzeit,izlog}. In the context of random
polynomials, there were many results on the universal character of
local correlations between close
zeros~\cite{bleherdi,ledoanetal,schiffmanzelditchbundles,shiffmanzelditch1}.
In this work, we focus on the global distribution of zeros.



The paper is organized as follows. In
Sections~\ref{subsecinvariant}--\ref{subsectheta}, we state our results
for a number of concrete ensembles of random analytic functions. These
results are special cases of the general Theorem~\ref{theogeneral}
whose statement, due to its technicality, is postponed to
Section~\ref{subsecgeneraltheo}. Proofs are given in
Sections~\ref{secproofsspecialcases} and~\ref{secproofgeneral}.

\subsection{Notation}\label{subsecnotation}
Let $\mathbb{D}_r=\{z\in\mathbb{C}\dvtx  |z|<r\}$ be the open disk with
radius \mbox{$r>0$} centered at the origin. Let
$\mathbb{D}=\mathbb{D}_1$ be the unit disk. Put
$\mathbb{D}_{\infty}=\mathbb{C}$. Denote by $\lambda$ the Lebesgue
measure on $\mathbb{C}$. A Borel measure $\mu$ on a locally compact
metric space $X$ is called \textit{locally finite} (l.f.) if
$\mu(A)<\infty$ for every compact set
$A\subset X$. 
A sequence $\mu_n$ of l.f. measures on $X$ converges \textit{vaguely}
to a l.f. measure $\mu$ if for every continuous, compactly supported
function $\varphi\dvtx X\to\mathbb{R}$,
%
\begin{equation}
\label{eqvagueconvmeas} 
\lim_{n\to\infty} \int
_{X} \varphi(z) \mu_n(dz) = \int
_{X} \varphi(z) \mu(dz). 
\end{equation}
If $\mu_n$ and $\mu$ are probability measures, the vague convergence is
equivalent to the more familiar weak convergence for
which~(\ref{eqvagueconvmeas}) is required to hold for all continuous,
bounded functions $\varphi$; see Lemma~4.20 in~\cite{kallenbergbook}.
Let $\mathcal{M}(X)$ be the space of all l.f. measures on $X$ endowed
with the vague topology. Note that $\mathcal{M}(X)$ is a Polish space;
see Theorem~A2.3 in~\cite{kallenbergbook}. A \textit{random measure} on
$X$ is a random element defined on some probability space
$(\Omega,\mathcal {F},\mathbb{P})$ and taking values in
$\mathcal{M}(X)$. The a.s. convergence and convergence in probability
of random measures are defined as the convergence of the corresponding
$\mathcal{M}(X)$-valued random elements. An equivalent definition: a
sequence of random measures $\mu_n$ converges to a random measure $\mu$
in probability (resp., a.s.), if~(\ref{eqvagueconvmeas}) holds in
probability (resp., a.s.) for every continuous, compactly supported
function $\varphi\dvtx X\to \mathbb{R}$.

%
%

\section{Statement of results}\label{secstatementofresults}
\subsection{The three invariant ensembles}\label{subsecinvariant}
Let $\xi_0,\xi_1,\ldots$ be i.i.d. random variables. Unless stated
otherwise, they take values in $\mathbb{C}$, are nondegenerate, and
satisfy the condition $\mathbb{E}\log(1+|\xi_0|)<\infty$. 
Fix a parameter $\alpha>0$. We start by considering the following three
ensembles of random analytic functions (see, e.g.,
\cite{sodintsirelson,peresbook}):
\[
\mathbf{F}_n(z) = \cases{ \displaystyle\sum
_{k=0}^n \xi_k \biggl(\frac{n(n-1)\cdots
(n-k+1)}{k!}
\biggr)^{\alpha} z^k &\quad(elliptic, $n\in\mathbb{N}$, $z\in
\mathbb{C}$), \vspace*{6pt}
\cr
\displaystyle\sum_{k=0}^{\infty}
\xi_k \biggl(\frac{n^k}{k!} \biggr)^{\alpha}
z^k &\quad(flat, $n>0$, $z\in\mathbb{C}$), \vspace*{6pt}
\cr
\displaystyle\sum_{k=0}^{\infty}
\xi_k \biggl( \frac{n(n+1)\cdots
(n+k-1)}{k!} \biggr)^{\alpha}
z^k &\quad(hyperbolic, $n>0$, $z\in\mathbb{D}$).}
\]
%
Note that in the elliptic case $\mathbf{F}_n$ is a random polynomial of
degree $n$, in the flat case it is a random entire function, whereas in
the hyperbolic case it is a random analytic function defined on the
unit disk $\mathbb{D}$. The a.s. convergence of the series in the
latter two
cases follows from Lemma~\ref{lemlogmoment} below. In the particular
case when $\alpha=1/2$ and $\xi_k$ are complex standard Gaussian with
density $z\mapsto\pi^{-1} \exp\{-|z|^2\}$ on $\mathbb{C}$, the zero
sets of
these analytic functions possess remarkable invariance properties
relating them to the three geometries of constant curvature;
see~\cite{sodintsirelson,peresbook}. In this special case, the expected
number of zeros of $\mathbf{F}_n$ in a Borel set $B$ can be computed
exactly~\cite{sodintsirelson,peresbook}:
%
\[
\mathbb{E} \bigl[\mu_{\mathbf{F}_n}(B) \bigr] = \cases{ \displaystyle\frac n{
\pi } \int_B \bigl(1+ |z|^{2}
\bigr)^{-2} \lambda(dz) &\quad(elliptic case, $B\subset\mathbb{C}$),
\vspace*{6pt}
\cr
\displaystyle \frac n{\pi} \lambda(B) &\quad(flat case, $B
\subset\mathbb{C}$), \vspace*{6pt}
\cr
\displaystyle\frac n{\pi} \int
_B \bigl(1-|z|^{2} \bigr)^{-2}
\lambda(dz) &\quad(hyperbolic case, $B\subset\mathbb{D}$).}
\]
%

In the next theorem, we compute the asymptotic distribution of zeros of
$\mathbf{F}_n$ for more general $\xi_k$'s.

%
\begin{theorem}\label{theoinvariant}
Let $\xi_0,\xi_1,\ldots$ be nondegenerate i.i.d. random variables such
that $\mathbb{E}\log(1+|\xi_0|)<\infty$. As $n\to\infty$, the sequence
of random measures $\frac1n \mu_{\mathbf{F}_n}$ converges in
probability to the deterministic measure having a density
$\rho_{\alpha}$ with respect to the Lebesgue measure, where
\[
\rho_{\alpha}(z) = \cases{ \displaystyle\frac1{2\pi\alpha
}|z|^{(1/\alpha) - 2} \bigl(1+ |z|^{1/\alpha} \bigr)^{-2} &
\quad(elliptic case, $z\in \mathbb{C}$), \vspace*{6pt}
\cr
\displaystyle\frac1{2
\pi\alpha} |z|^{(1/\alpha)-2} &\quad (flat case, $z\in\mathbb{C}$), \vspace*{6pt}
\cr
\displaystyle\frac1 {2\pi\alpha} |z|^{(1/\alpha) - 2} \bigl(1-|z|^{1/\alpha}
\bigr)^{-2} & \quad(hyperbolic case, $z\in\mathbb{D}$).}
\]
%
\end{theorem}


\subsection{Littlewood--Offord random polynomials}\label{subsecLOpoly}
Next, we consider an ensemble of random polynomials which was
introduced by Littlewood and Offord
\cite{littlewoodofford1,littlewoodofford2}. It is related to the flat
model. First, we give some motivation. Let $\xi_0,\xi_1,\ldots$ be
nondegenerate i.i.d. random variables. Given a sequence
$w_0,w_1,\ldots\in\mathbb{C}\setminus\{0\}$ consider a random
polynomial $\mathbf{W}_n$ defined by
%
\begin{equation}
\label{eqdefQnLO} \mathbf{W}_n(z)=\sum_{k=0}^n
\xi_k w_k z^k.
\end{equation}
For $w_k=1$, we recover the Kac polynomials, for which the zeros
concentrate near the unit circle. The next result shows that the
structure of the zeros does not differ essentially from the Kac case if
the sequence $w_k$ grows or decays not too fast.

%
\begin{theorem}\label{theokacgeneralized}
Let $\xi_0,\xi_1,\ldots$ be nondegenerate i.i.d. random variables such
that $\mathbb{E}\log(1+|\xi_0|)<\infty$. If \ $\lim_{k\to\infty}
\frac1k \log |w_k| = w$ for some constant $w\in\mathbb{R}$, then the
sequence of random measures $\frac1n \mu_{\mathbf{W}_n}$ converges in
probability to the uniform probability distribution on the circle of
radius $e^{-w}$ centered at the origin.
\end{theorem}
%

We would like to construct examples where there is no concentration
near a circle.
Let us make the following assumption on the sequence $w_k$:
%
\begin{equation}
\label{eqfkasympt} \log|w_k| = - \alpha(k\log k - k) - \beta k
+o(k), \qquad k\to\infty,
\end{equation}
where $\alpha>0$ and $\beta\in\mathbb{R}$ are parameters.
Particular cases are
polynomials of the form
\begin{eqnarray*}
\mathbf{W}_n^{(1)}(z)&=&\sum_{k=0}^n \frac{\xi_k}{(k!)^{\alpha}}
z^k,
\\
\mathbf{W}_n^{(2)}(z)&=& \sum_{k=0}^n \frac{\xi_k}{k^{\alpha
k}}z^k,
\\
\mathbf{W}_n^{(3)}(z)&=&\sum _{k=0}^n
\frac{\xi_k}{\Gamma(\alpha k+1)} z^k.
\end{eqnarray*}
The family $\mathbf{W}_n^{(1)}$ has been studied by Littlewood and Offord
\cite{littlewoodofford1,littlewoodofford2} in one of the earliest works
on random polynomials. They were interested in the number of real
zeros. 
In the next theorem, we describe the limiting distribution of complex
zeros of $\mathbf{W}_n$. Let $\mu_n$ be the measure counting the
points of
the form $e^{-\beta} n^{-\alpha}z$, where $z$~is a zero of $\mathbf{W}_n$.
That is, for every Borel set $B\subset\mathbb{C}$,
%
\begin{equation}
\label{eqmuwwwnnormalized} \mu_n(B)=\mu_{\mathbf{W}_n}
\bigl(e^{\beta} n^{\alpha}B \bigr).
\end{equation}
%

%
\begin{theorem}\label{theoLOpoly}
Let $\xi_0,\xi_1,\ldots$ be nondegenerate i.i.d. random variables such
that $\mathbb{E}\log(1+|\xi_0|)<\infty$. Let $w_0,w_1,\ldots$ be a
complex sequence satisfying~(\ref{eqfkasympt}). With probability $1$,
the sequence of random measures $\frac1n \mu_n$ converges to the
deterministic probability measure having the density
%
\begin{equation}
\label{eqLOpolydensity} z\mapsto\frac1{2\pi\alpha} |z|^{(1/\alpha)-2}
\mathbh{1}_{z\in\mathbb{D}}
\end{equation}
with respect to the Lebesgue measure on $\mathbb{C}$.
\end{theorem}

%
For the so-called \textit{Weyl random polynomials}
$\mathbf{W}_n(z)=\sum_{k=0}^n \xi_k \frac{z^k}{\sqrt{k!}}$ having
$\alpha=1/2$ and $\beta=0$, the limiting distribution is uniform on
$\mathbb{D}$; see Figure~\ref{figweyluniversality}. This result can
be seen
as an analogue of the famous circular law for the distribution of
eigenvalues of non-Hermitian random matrices with i.i.d.
entries~\cite{taovu,bordenavechafai}. Forrester and Honner
\cite{forresterhonner} stated the circular law for Weyl polynomials and
discussed differences and similarities between the matrix and the
polynomial cases; see also~\cite{krishnapurvirag}.


Under a minor additional assumption on the coefficients $w_k$ we can
prove that the logarithmic moment condition is not only sufficient, but
also necessary for the a.s. convergence of the empirical distribution
of zeros. It is easy to check that the additional assumption is
satisfied for $\mathbf{W}_n=\mathbf{W}_n^{(i)}$ with $i=1,2,3$.

%
\begin{theorem}\label{theoLOpolyconverse}
Let $\xi_0,\xi_1,\ldots$ be nondegenerate i.i.d. random variables. Let
$w_0,w_1,\ldots$ be a complex sequence satisfying~(\ref{eqfkasympt})
and such that for some $C>0$,
%
\begin{equation}
\label{eqLOadditionalassumpt} |w_{n-k}/w_n|< C
e^{\beta k} n^{\alpha k}\qquad\mbox{ for all } n\in\mathbb{N}, k\leq n.
\end{equation}
Let $\mu_n$ be as in~(\ref{eqmuwwwnnormalized}). Then, the following
are equivalent:
\begin{longlist}[(2)]
\item[(1)] With probability $1$, the sequence of random measures $\frac
1n\mu_n$ converges to the probability measure with
density~(\ref{eqLOpolydensity}).

\item[(2)] $\mathbb{E}\log(1+|\xi_0|)<\infty$.
\end{longlist}
\end{theorem}

\begin{figure}

\includegraphics{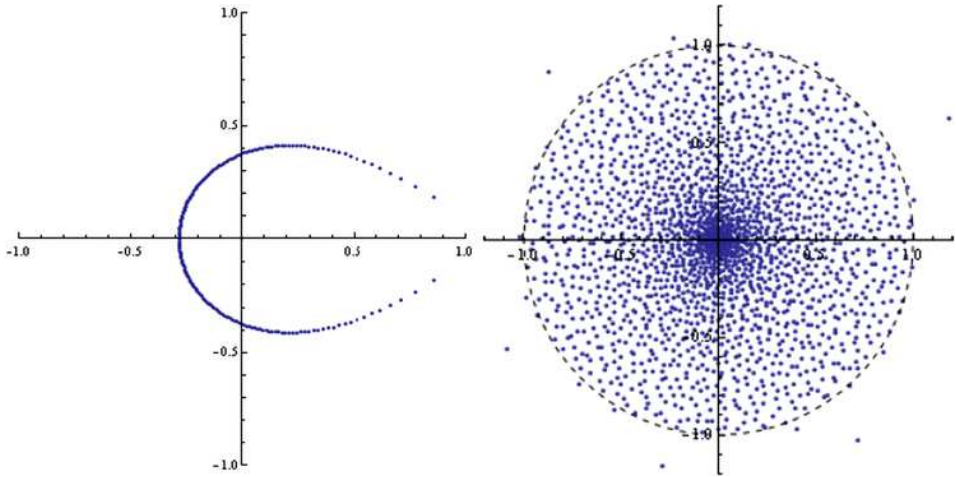}

\caption{Left: Zeros of the Szeg\H{o} polynomial $s_n(z)=\sum_{k=0}^n
\frac{z^k}{k!}$ of degree $n=200$. Right: Zeros of the
Littlewood--Offord random polynomial $\mathbf{W}_n(z)=\sum_{k=0}^n
\xi_k\frac{z^k}{k!}$ of degree $n=2000$ with complex normal coefficients. In both
cases, the zeros were divided by $n$.}\label{figszegoe}
\end{figure}

It should be stressed that in all our results we assume that the random
variables $\xi_k$ are nondegenerate (i.e., not a.s. constant). To see
that this assumption is essential, consider the deterministic
polynomials
%
\begin{equation}
\label{eqszegoepoly} s_n(z)=\sum_{k=0}^n
\frac{z^k}{k!}.
\end{equation}
A classical result of Szeg{\H{o}} \cite{szegoe} states that the zeros
of $s_n(nz)$ cluster asymptotically (as $n\to\infty$) along the curve
$\{|ze^{1-z}|=1\}\cap\mathbb{D}$; see Figure~\ref{figszegoe}~(left).
This behavior is manifestly different from the distribution with
density $1/(2\pi|z|)$ on~$\mathbb{D}$ we have obtained in
Theorem~\ref{theoLOpoly} for the same polynomial with randomized
coefficients; see Figure~\ref{figszegoe}~(right).
%
%

\subsection{Littlewood--Offord random entire function}\label{subsecLOentire}
Next we discuss a random entire function which also was introduced by
Littlewood and Offord
\cite{littlewoodoffordentire1,littlewoodoffordentire2err}. Their aim
was to describe the properties of a ``typical'' entire function of a
given
order~$1/\alpha$. 
Given a complex sequence $w_0,w_1,\ldots$ satisfying~(\ref{eqfkasympt})
consider a random entire function
%
\begin{equation}
\label{eqdefLOentire} \mathbf{W}(z)=\sum_{k=0}^{\infty}
\xi_k w_k z^k.
\end{equation}
Examples are given by
\begin{eqnarray*}
\mathbf{W}^{(1)}(z)&=&\sum_{k=0}^{\infty}
\frac{\xi_k}{(k!)^{\alpha}} z^k,
\\
\mathbf{W}^{(2)}(z)&=&\sum_{k=0}^{\infty} \frac{\xi_k}{k^{\alpha k}}z^k,
\\
\mathbf{W}^{(3)}(z)&=& \sum_{k=0}^{\infty}\frac{\xi_k}{\Gamma(\alpha k+1)} z^k.
\end{eqnarray*}
The first function is essentially the flat model considered above,
namely $\mathbf{W}^{(1)}(n^{\alpha}z)=\mathbf{F}_n(z)$. For $\alpha
=1$, it is a
randomized version of the Taylor series for the exponential. The last
function is a randomized version of the Mittag--Leffler function.
Our aim is to describe the density of zeros of $\mathbf{W}$ on the global
scale.
Let $\mu_n$ be the measure counting the points of the form $e^{-\beta}
n^{-\alpha}z$, where $z$ is a zero of~$\mathbf{W}$. That is, for
every Borel
set $B\subset\mathbb{C}$,
%
\begin{equation}
\label{eqmuwwwnormalized} \mu_n(B)= \mu_{\mathbf{W}}
\bigl(e^{\beta} n^{\alpha}B \bigr).
\end{equation}
We have the following strengthening of the flat case of
Theorem~\ref{theoinvariant}.

%
\begin{theorem}\label{theoLOentire}
Let $\xi_0,\xi_1,\ldots$ be nondegenerate i.i.d. random variables such
that $\mathbb{E}\log(1+|\xi_0|)<\infty$. Let $w_0,w_1,\ldots$ be a
complex sequence satisfying~(\ref{eqfkasympt}). With probability $1$,
the random measure $\frac1n \mu_{n}$ converges to the deterministic
measure having the density
%
\begin{equation}
\label{eqLOentiredensity} z\mapsto\frac1{2\pi\alpha} |z|^{(1/\alpha)-2}
\end{equation}
with respect to the Lebesgue measure on $\mathbb{C}$.
\end{theorem}
As a corollary, we obtain a law of large numbers for the number of
zeros of $\mathbf{W}$.

%
\begin{corollary}\label{corLOentireLLN}
Let $N(r)=\mu_{\mathbf{W}}(\mathbb{D}_r)$ be the number of zeros of
$\mathbf{W} $ in the disk $\mathbb{D}_r$. Under the assumptions of
Theorem~\ref{theoLOentire},
\[
N(r) = e^{-\beta/\alpha} r^{1/\alpha} \bigl(1+o(1)
\bigr)\qquad\mbox{a.s. as } r\to\infty.
\]
\end{corollary}

In the case $\alpha=1/2$ the limiting measure in
Theorem~\ref{theoLOentire} has constant density~$1/\pi$. The difference
between the limiting densities in Theorems
\ref{theoLOpoly}~and~\ref{theoLOentire} is that in the latter case
there is no restriction to the unit disk.
It has been pointed out by the unknown referee that in the special case
of the Bernoulli-distributed $\xi_k$'s Theorem~\ref{theoLOentire} can
be deduced from the results of Littlewood and Offord
\cite{littlewoodoffordentire1,littlewoodoffordentire2err} using the
Levin--Pfluger theory (\cite{levin}, Chapter~3). Our proof is simpler
than the proof of Littlewood and Offord
\cite{littlewoodoffordentire1,littlewoodoffordentire2err}. For a
related work, see also~\cite{offord1,offord2}.



Let us again stress the importance of the nondegeneracy assumption. The
exponential function $e^z$ has no complex zeros, whereas the zeros of
its randomized version $\sum_{k=0}^{\infty}\xi_k \frac{z^k}{k!}$ have
the global-scale density $1/(2\pi|z|)$ on $\mathbb{C}$. For the
absolute values of the zeros, the limiting density is constant and
equal to $1$ on $(0,\infty)$.

\subsection{Randomized theta function}\label{subsectheta}
Given a parameter $\alpha\in(0,1)\cup(1,\infty)$ we consider a random
analytic function
\[
\mathbf{H}_n(z) = \cases{ \displaystyle\sum
_{k=0}^{\infty} \xi_k e^{n^{1-\alpha} k^{\alpha}}
z^k &\quad(case $\alpha<1$, $z\in \mathbb{D}$), \vspace*{6pt}
\cr
\displaystyle\sum_{k=0}^{\infty}
\xi_k e^{-n^{1-\alpha} k^{\alpha}} z^k &\quad(case $\alpha>1$, $z\in
\mathbb{C}$).}
\]

%
\begin{theorem}\label{theostable}
Let $\xi_0,\xi_1,\ldots$ be nondegenerate i.i.d. random variables such
that $\mathbb{E}\log(1+|\xi_0|)<\infty$. As $n\to\infty$, the sequence
of random measures $\frac1n \mu_{\mathbf{H}_n}$ converges in
probability to the deterministic measure having the density
\[
z \mapsto\frac{1}{2\pi\alpha|1-\alpha|} \frac1 {|z|^2}
\biggl|\frac{\log|z|}{\alpha} \biggr|^{(2-\alpha)/(\alpha-1)}
\]
with respect to the Lebesgue measure on $\mathbb{C}$. The density is
restricted to $\mathbb{D}$ in the case $\alpha<1$ and to
$\mathbb{C}\setminus \mathbb{D}$ in the case $\alpha>1$.
\end{theorem}

As the parameter $\alpha$ crosses the value $1$, the zeros of
$\mathbf{H}_n$ jump from the unit disk~$\mathbb{D}$ to its complement
$\mathbb {C}\setminus\mathbb{D}$. Note that the case $\alpha=1$
corresponds formally to Kac polynomials for which the zeros are on the
boundary of $\mathbb{D}$. The special case $\alpha=2$ corresponds to
the randomized theta function
%
\begin{equation}
\label{eqdeftheta} \mathbf{H}_n(z)=\sum
_{k=0}^{\infty} \xi_k e^{-k^2/n}
z^k.
\end{equation}
The limiting distribution of zeros has the density $\frac1{4\pi|z|^2}$
on $\mathbb{C}\setminus\mathbb{D}$. One can also take the sum
in~(\ref{eqdeftheta}) over $k\in\mathbb{Z}$ in which case the zeros
fill the whole complex plane with the same density.

A similar model, namely the polynomials $\mathbf{Q}_n(z)=\sum_{k=0}^n
\xi_k e^{-k^{\alpha}} z^k$, where \mbox{$\alpha>1$}, has been considered by
Schehr and Majumdar \cite{schehrmajumdar}. Assuming that $\xi_k$ are
real-valued they showed that almost all zeros of $\mathbf{Q}_n$ become
real if $\alpha>2$. In our model, the distribution of the arguments of
the zeros remains uniform for every $\alpha$.


\subsection{The general result}\label{subsecgeneraltheo}
We are going to state a theorem which contains all examples considered
above as special cases. Let $\xi_0,\xi_1,\ldots$ be nondegenerate
i.i.d. complex-valued random variables such that $\mathbb{E}
\log(1+|\xi_0|)<\infty$.
Consider a random Taylor series
%
\begin{equation}
\label{eqdefGGG} \mathbf{G}_n(z)=\sum_{k=0}^{\infty}
\xi_k f_{k,n} z^k,
\end{equation}
where $f_{k,n}\in\mathbb{C}$ are deterministic coefficients.
Essentially, we
will assume that for some function $u(t)$ the coefficients $f_{k,n}$
satisfy
\[
|f_{k,n}| = e^{-n u(k/n)+ o(n)},\qquad n\to\infty.
\]
Here is a precise statement. We assume that there is a function $f\dvtx
[0,\infty)\to[0,\infty)$ and a number $T_0\in(0,\infty]$ such that
\begin{longlist}[(A3)]
\item[(A1)] $f(t)>0$ for $t<T_0$ and $f(t)=0$ for $t>T_0$.

\item[(A2)] $f$ is continuous on $[0,T_0)$, and, in the case
$T_0<+\infty$, left continuous at~$T_0$.

\item[(A3)]
$ \lim_{n\to\infty} \sup_{k\in[0,An]} \vert
|f_{k,n}|^{1/n}-f(\frac kn)\vert=0
$ for every $A>0$.

\item[(A4)] $R_0:=\liminf_{t\to\infty} f(t)^{-1/t}\in(0,\infty]$,
    $\liminf_{k\to\infty} |f_{k,n}|^{-1/k} \geq R_0$ for every fixed
    $n\in\mathbb{N}$ and additionally, $\liminf_{n,k/n\to \infty}
    |f_{k,n}|^{-1/k} \geq R_0$.
\end{longlist}
%
It will be shown later that condition~(A4) ensures that the
series~(\ref{eqdefGGG}) defining $\mathbf{G}_n$ converges with
probability $1$ on the disk $\mathbb{D}_{R_0}$. Let $I\dvtx
\mathbb{R}\to\mathbb {R}\cup\{+\infty\} $ be the Legendre--Fenchel
transform of the function $u(t)=-\log f(t)$, where $\log0=-\infty$.
That is,
%
\begin{equation}
\label{eqdefI} I(s)=\sup_{t\geq0} \bigl(st-u(t) \bigr)=\sup
_{t\geq0} \bigl(st+\log f(t) \bigr). 
\end{equation}
Note that $I$ is a convex function, $I(s)$ is finite for $s<\log R_0$
and $I(s)=+\infty$ for $s>\log R_0$. Recall that $\mu_{\mathbf
{G}_n}$ is the
measure assigning to each zero of $\mathbf{G}_n$ a weight equal to its
multiplicity.

%
\begin{theorem}\label{theogeneral}
Under the above assumptions, the sequence of random measures $\frac1n
\mu_{\mathbf{G}_n}$ converges in probability to some deterministic locally
finite measure $\mu$ on the disk $\mathbb{D}_{R_0}$. The measure $\mu
$ is
rotationally invariant and is characterized by
%
\begin{equation}
\label{eqtheogeneral} 
\mu(\mathbb{D}_r)
=I'( \log r),\qquad r\in(0,R_0). 
\end{equation}
%
\end{theorem}

By convention, $I'$ is the left derivative of $I$. Since $I$ is convex,
the left derivative exists everywhere on $(-\infty, \log R_0)$ and is a
nondecreasing, left-continuous function.
Since the supremum in~(\ref{eqdefI}) is taken over $t\geq0$, we have
$\lim_{s\to-\infty}I'(s)=0$. Hence, $\mu$ has no atom at zero. If $I'$
is absolutely continuous on some interval $(\log r_1,\log r_2)$, then
the density of $\mu$ on the annulus $r_1<|z|<r_2$ with respect to the
Lebesgue measure on $\mathbb{C}$ is
%
\begin{equation}
\label{eqrhoformula} \rho(z)=\frac{I''(\log|z|)} {2\pi|z|^2}.
\end{equation}
%

It is possible to give a characterization of the measure $\mu$ without
referring to the Legendre--Fenchel transform. The radial part of $\mu$
is a measure $\bar\mu$ on $(0,\infty)$ defined by $\bar
\mu((0,r))=\mu(\mathbb{D}_r)$. Suppose first that $u$ is convex on
$(0,T_0)$ (which is the case in all our examples). Then, $\bar\mu$ is
the image of the Lebesgue measure on $(0,\infty)$ under the mapping
$t\mapsto e^{u'(t)}$, where $u'$ is the left derivative of $u$. This
follows from the fact that $(u')^{\leftarrow}=I'$ and
$(I')^{\leftarrow}=u'$ by the Legendre--Fenchel duality, where
$\varphi^{\leftarrow}(t)=\inf\{s\in \mathbb{R}\dvtx  \varphi(s)\geq
t\}$ is the generalized left-continuous inverse of a nondecreasing
function $\varphi$. In particular, the support of $\mu$ is contained in
the annulus
\[
\bigl\{e^{\lim_{t\downarrow0}u'(t)}\leq|z|\leq e^{\lim_{t\uparrow
T_0} u'(t)} \bigr\}
\]
and is equal to this annulus if $u'$ has no jumps.
In general, any jump of $u'$ (or, by duality, any constancy interval of
$I'$) corresponds to a missing annulus in the support of $\mu$. Also,
any jump of $I'$ (or, by duality, any constancy interval of $u'$)
corresponds to a circle with positive $\mu$-measure. More precisely, if
$I'$ has a jump at $s$ (or, by duality, $u'$ takes the value $s$ on an
interval of positive length), then $\mu$ assigns a positive weight
(equal to the size of the jump) to the circle of radius $e^s$ centered
at the origin. In the case when $u$ is nonconvex, we can apply the same
considerations after replacing $u$ by its convex hull.

One may ask what measures $\mu$ may appear as limits in
Theorem~\ref{theogeneral}. Clearly, $\mu$ has to be rotationally
invariant, with no atom at $0$. The next theorem shows that there are
no further essential restrictions.

%
\begin{theorem}\label{theogeneralconverse}
Let $\mu$ be a rotationally invariant measure on $\mathbb{C}$ such that
\begin{longlist}[(2)]
\item[(1)] $\mu(\mathbb{C}\setminus\mathbb{D}_{R_0})=0$, where
$R_0:=\sup\{r>0\dvtx \mu(\mathbb{D}_r)<\infty\}\in(0,\infty]$.

\item[(2)] $\int_{0}^{R}\mu(\mathbb{D}_r)r^{-1} \,dr<\infty$ for some
(hence, every) $R<R_0$.
\end{longlist}
Then, there is a random Taylor series $\mathbf{G}_n$ of the
form~(\ref{eqdefGGG}) with convergence radius a.s. $R_0$ such that
$\frac1 n\mu_{\mathbf{G}_n}$ converges in probability to $\mu$ on the
disk $\mathbb{D}_{R_0}$.
\end{theorem}

%
%
\begin{example}\label{exthreecircles}
Consider a random polynomial
%
\begin{equation}
\label{eqdefthreecircles} \mathbf{G}_n(z)=\sum
_{k=0}^{n} \xi_k z^k +
2^n\sum_{k=n+1}^{2n}
\xi_k \biggl(\frac{z}{2} \biggr)^k + \biggl(\frac9
2 \biggr)^n \sum_{k=2n+1}^{3n}
\xi_k \biggl(\frac
{z}{3} \biggr)^k.
\end{equation}
We can apply Theorem~\ref{theogeneral} with
\begin{eqnarray*}
u(t) &=& \cases{ 0, &\quad$t\in[0,1]$, \vspace*{2pt}
\cr
(\log2) (t-1), &\quad$t
\in[1,2]$, \vspace*{2pt}
\cr
(\log3) t -\log\frac92, &\quad$t\in[2,3]$,
\vspace*{2pt}
\cr
{+}\infty, &\quad$t\geq3$,}
\\
I(s) &=& \cases{ 0, &\quad$s
\leq0$, \vspace*{2pt}
\cr
s, &\quad$s\in[0,\log2]$, \vspace*{2pt}
\cr
2s -\log2, &
\quad$s\in[\log2,\log3]$, \vspace*{2pt}
\cr
3s -\log6, &\quad$t\geq\log3$.}
\end{eqnarray*}
The function $u'$ has three constancy intervals of length $1$ where it
takes values $0,\log2,\log3$. Dually, the function $I'$ has three
jumps of size $1$ at $0,\log2, \log3$ and is locally constant outside
these points. It follows that the limiting distribution of the zeros of
$\mathbf{G}_n$ is the sum of uniform probability distributions on three
concentric circles with radii $1,2,3$.
\end{example}

%
\begin{remark}
Suppose that $\mathbf{G}_n$ satisfies the assumptions of
Theorem~\ref{theogeneral}. Then, so does the derivative $\mathbf
{G}_n'$ (and, moreover, $f$ is the same in both cases). Thus, the
derivative of any fixed order of $\mathbf{G}_n$ has the same limiting
distribution of zeros as $\mathbf{G}_n$. Similarly, for every complex
sequence $c_n$ such that $\limsup_{n\to\infty} \frac1n \log|c_n|\leq
f(0)$, the function $\mathbf{G}_n(z)-c_n$ satisfies the assumptions.
Hence, the limiting distribution of the solutions of the equation
$\mathbf{G}_n(z)=c_n$ is the same as for the zeros of~$\mathbf{G}_n$.
\end{remark}

\section{Proofs: Special cases}\label{secproofsspecialcases}
We are going to prove the results of Section~\ref{secintro}. We will
verify the assumptions of Section~\ref{subsecgeneraltheo} and apply
Theorem~\ref{theogeneral}. Recall the notation $u(t)=-\log f(t)$.

\begin{pf*}{Proof of Theorem~\ref{theokacgeneralized}}
We can assume that $w=0$ since otherwise we can consider the polynomial
$\mathbf{W}_n(e^{-w}z)$. It follows from $\lim_{k\to\infty}\frac1k \log
|w_k|=0$ that assumptions~(A1)--(A4) of Section~\ref{subsecgeneraltheo}
are fulfilled with $T_0=1$, $R_0=+\infty$ and
\[
f(t)= \cases{ 1, &\quad$t\in[0,1]$, \vspace*{2pt}
\cr
0, &\quad$t>1$,} \qquad
u(t)= \cases{ 0, &\quad$t\in[0,1]$, \vspace*{2pt}
\cr
{+}\infty, &\quad$t>1$.}
\]
%
The Legendre--Fenchel transform of $u$ is given by $I(r)=\max(0, r)$.
It follows from~(\ref{eqtheogeneral}) that $\mu$ is the uniform
probability measure on $\mathbb{T}$.
\end{pf*}

%
\begin{remark}
Under a slightly more restrictive assumption $\mathbb{E}\log|\xi
_0|<\infty$, Theorem~\ref{theokacgeneralized} can be deduced from the
result of Hughes and Nikeghbali \cite{hughesnikeghbali} (which is
partially based on the Erd\H{o}s--Turan inequality). This method,
however, requires a subexponential growth of the coefficients and
therefore fails in all other examples we consider here.
\end{remark}

\begin{pf*}{Proof of Theorem~\ref{theoinvariant}}
By the Stirling formula, $\log n!=n\log n-n+o(n)$ as $n\to\infty$. It
follows that assumption~(A3) holds with
\[
u(t)= \cases{ \alpha \bigl(t\log t+ (1-t)\log(1-t) \bigr), &\quad $0\leq t
\leq1$, elliptic case, \vspace*{2pt}
\cr
\alpha(t\log t-t), &\quad $t\geq0$, flat
case, \vspace*{2pt}
\cr
\alpha \bigl(t\log t- (1+t)\log(1+t) \bigr), &\quad $t
\geq0$, hyperbolic case.}
\]
In the elliptic case, $u(t)=+\infty$ for $t>1$. The Legendre--Fenchel
transform of $u$ is given by
\[
I(s)= \cases{ \alpha\log \bigl(1+e^{s/\alpha} \bigr), &\quad $s\in
\mathbb{R}$, elliptic case, \vspace*{2pt}
\cr
\alpha e^{s/\alpha}, &\quad $s
\in\mathbb{R}$, flat case, \vspace*{2pt}
\cr
{-}\alpha\log \bigl(1-e^{s/\alpha}
\bigr), & \quad $s<0$, hyperbolic case.}
\]
In the hyperbolic case, $I(s)=+\infty$ for $s\geq0$. We have $R_0=1$ in
the hyperbolic case and $R_0=+\infty$ in the remaining two cases. The
proof is completed by applying Theorem~\ref{theogeneral}.
\end{pf*}

\begin{pf*}{Proof of Theorem~\ref{theoLOpoly}}
We are going to apply Theorem~\ref{theogeneral} to the polynomial
$\mathbf{G}_n(z)=\mathbf{W}_n(e^{\beta} n^{\alpha} z)$. We have
$f_{k,n}=e^{\beta k
+ \alpha k \log n} w_{k}$, for $0\leq k\leq n$.
Equation~(\ref{eqfkasympt}) implies that assumption~(A3) is satisfied
with
\[
u(t) = \cases{ \alpha(t\log t-t), &\quad$t\in[0,1]$, \vspace*{2pt}
\cr
{+}
\infty, &\quad$t>1$.}
\]
%
The Legendre--Fenchel transform of $u$ is given by
\[
I(s)= \cases{ \alpha e^{s/\alpha}, &\quad$s\leq0$, \vspace*{2pt}
\cr
\alpha+
s, &\quad$s\geq0$.}
\]
Applying Theorem~\ref{theogeneral}, we obtain that $\frac1n
\mu_{\mathbf{G}_n}$ converges in probability to the required limit. A.s.
convergence will be demonstrated in Section~\ref{subsecproofLOpolyas}
below.
\end{pf*}

\begin{pf*}{Proof of Theorem~\ref{theoLOentire}}
We apply Theorem~\ref{theogeneral} to $\mathbf{G}_n(z)=\mathbf
{W}(e^{\beta} n^{\alpha} z)$. We have $u(t)=\alpha(t\log t-t)$ for all
$t\geq0$. Hence, $I(s)=\alpha e^{s/\alpha}$ for all $s\in\mathbb{R}$.
We can apply Theorem~\ref{theogeneral} to prove convergence in
probability. A.s. convergence will be demonstrated in
Section~\ref{subsecproofLOentireas} below.
\end{pf*}

\begin{pf*}{Proof of Theorem~\ref{theostable}}
Put $\sigma=+1$ in the case $\alpha>1$ and $\sigma=-1$ in the case
$\alpha<1$.
We have $u(t)=\sigma t^{\alpha}$ for $t\geq0$. It follows that
\[
I(r) = \cases{ \displaystyle\sigma(\alpha-1) \biggl(\frac{\sigma
r}{\alpha}
\biggr)^{\alpha/(\alpha-1)}, &\quad$\sigma r\geq0$, \vspace*{2pt}
\cr
{+}\infty, &
\quad$ \sigma r<0$.}
\]
We can apply Theorem~\ref{theogeneral}.
\end{pf*}

\section{Proofs: General results}\label{secproofgeneral}
\subsection{Method of proof of Theorem~\texorpdfstring{\protect\ref{theogeneral}}{2.8}}
We use the notation and the assumptions of
Section~\ref{subsecgeneraltheo}. We denote the probability space on
which the random variables $\xi_0,\xi_1,\ldots$ are defined by
$(\Omega, \mathcal{F}, \mathbb{P})$. We will write $\mu_n=\mu
_{\mathbf{G}_n}$ for the measure counting the zeros of $\mathbf{G}_n$.
To stress the randomness of the object under consideration we will
sometimes write $\mathbf{G}_n(z; \omega)$ and $\mu_n(\omega)$ instead
of $\mathbf{G}_n(z)$ and $\mu_{n}$. Here, $\omega\in\Omega$. The
starting point of the proof of Theorem~\ref{theogeneral} is the formula
%
\begin{equation}
\label{eqmuGGGisLaplacelogGGG} \mu_{n}(\omega)=\frac1 {2\pi} \Delta
\log \bigl|\mathbf{G}_n(z; \omega) \bigr|
\end{equation}
for every fixed $\omega\in\Omega$ for which $\mathbf{G}_n(z;\omega)$
does not vanish identically. Here, $\Delta$~denotes the Laplace
operator in the complex $z$-plane. The Laplace operator should always
be understood as an operator acting on
$\mathcal{D}'(\mathbb{D}_{R_0})$, the space of generalized functions on
the disk $\mathbb{D}_{R_0}$; see, for example, Chapter~II
of~\cite{hoermanderbook1}. Equation~(\ref{eqmuGGGisLaplacelogGGG})
follows from the formula $\frac1 {2\pi}\Delta\log|z-z_0|=\delta(z_0)$,
for every $z_0\in\mathbb{C}$; see Example~4.1.10
in~\cite{hoermanderbook1}. First, we will compute the limiting
logarithmic potential in~(\ref{eqmuGGGisLaplacelogGGG}).

%
\begin{theorem}\label{theologasympt}
Under the assumptions of Section~\ref{subsecgeneraltheo}, for every
$z\in\mathbb{D}_{R_0}\setminus\{0\}$,
%
\begin{equation}
\label{eqlogasympt} p_n(z):=\frac1n \log \bigl| \mathbf{G}_n(z)
\bigr| \mathop{\longrightarrow}^{P}_{n\to\infty}I \bigl(\log|z| \bigr).
\end{equation}
\end{theorem}

%
We will prove Theorem~\ref{theologasympt} in
Sections~\ref{subseclogmomentcond}, \ref{subsecupperbound},
\ref{subseclowerbound} below. Theorem~\ref{theologasympt} follows from
equations~(\ref{equpperboundG}) and~(\ref{eqprooflogpartfunc2}) below.
Moreover, it follows from~(\ref{equpperboundG}) that
$\limsup_{n\to\infty} p_n(z)\leq I(\log|z|)$ a.s. Unfortunately, we
were unable to prove that $\liminf_{n\to\infty} p_n(z)\geq I(\log|z|)$
a.s. Instead, we have the following slightly weaker statement.

%
\begin{proposition}\label{proplogasymptsubseq}
Let $l_1,l_2,\ldots$ be an increasing sequence of natural numbers such
that $l_k\geq k^3$ for all $k\in\mathbb{N}$. Under the assumptions of
Section~\ref{subsecgeneraltheo} we have, for every $z\in
\mathbb{D}_{R_0}\setminus\{0\}$,
%
\begin{equation}
\label{eqplktoI} p_{l_k}(z)=\frac1{l_k} \log \bigl|
\mathbf{G}_{l_k}(z) \bigr| \mathop {\longrightarrow}^{\mathrm{a.s.}}_{k\to\infty}
I \bigl(\log|z| \bigr).
\end{equation}
\end{proposition}
Proposition~\ref{proplogasymptsubseq} follows from
equations~(\ref{equpperboundG}) and~(\ref{eqprooflogpartfunc2}) by
noting that\break  \mbox{$\sum_{k=1}^{\infty} k^{-3/2}<\infty$} and applying the
Borel--Cantelli lemma. The next proposition allows us to pass from
convergence of potentials to convergence of measures. We will prove it
Section~\ref{subsecconvpotconvmeas}. Recall that $\mu_n$ counts the
zeros of~$\mathbf{G}_n$.

%
\begin{proposition}\label{lemconvpotconvmeas}
Let $l_1,l_2,\ldots$ be any increasing sequence of natural numbers.
Assume that for Lebesgue-a.e. $z\in\mathbb{D}_{R_0}$
equation~(\ref{eqplktoI}) holds. Then,
%
\begin{equation}
\label{eqlimitmulk} \frac{1}{l_k}\mu_{l_k} \mathop{
\longrightarrow}^{\mathrm
{a.s.}}_{k\to\infty}\frac{1}{2\pi} \Delta I \bigl(
\log|z| \bigr). 
\end{equation}
\end{proposition}

With these results, we are in position to prove
Theorem~\ref{theogeneral}. We need to show that $\frac1n \mu_{n}$
converges to $\mu$ in probability, as a sequence of
$\mathcal{M}(\mathbb{D}_{R_0})$-valued random variables. A sequence of
random variables with values in a metric space \mbox{converges} in probability
to some limit if and only if every subsequence of these random
variables contains a subsubsequence which converges a.s. to the same
limit; see, for example, Lemma~3.2 in~\cite{kallenbergbook}. Let a
subsequence $\frac{1}{n_1}\mu_{n_1}, \frac{1}{n_2}\mu_{n_2},\ldots,$
where $n_1<n_2<\cdots$, be given. Write $l_k=n_{k^3}$, so that
$\{l_k\}$ is a subsequence of $\{n_k\}$ and $l_k\geq k^3$. It follows
from Propositions
\ref{proplogasymptsubseq}~and~\ref{lemconvpotconvmeas}
that~(\ref{eqlimitmulk}) holds. So, the random measure $\frac1n \mu_n$
converges in probability to $\frac{1}{2\pi} \Delta I(\log|z|)$. It
remains to observe that the generalized function $\frac{1}{2\pi} \Delta
I(\log|z|)$ is equal to the measure $\mu$ given
in~(\ref{eqtheogeneral}). This follows from the fact that the radial
part of $\Delta$ in polar coordinates is given by $\frac1r \frac d {dr}
r \frac d {dr}$. This gives the desired result.

\subsection{The logarithmic moment condition} \label{subseclogmomentcond}
The next well-known lemma states that i.i.d. random variables grow
subexponentially with probability $1$ if and only if their logarithmic
moment is finite.

%
\begin{lemma}\label{lemlogmoment}
Let $\xi_0,\xi_1,\ldots$ be i.i.d. random variables. Fix
$\varepsilon>0$.
Then,
%
\begin{equation}
\label{eqdefmaxxik} S:=\sup_{k=0,1,\ldots} \frac{|\xi_k|}{e^{\varepsilon k}} < +
\infty\mbox{ a.s.}\quad \Longleftrightarrow\quad\mathbb{E}\log \bigl(1+|
\xi_0| \bigr)<\infty.
\end{equation}
\end{lemma}

\begin{pf}
For every nonnegative random variable $X$ we have
\[
\sum_{k=1}^{\infty} \mathbb{P}[X\geq k] \leq
\mathbb{E}X\leq \sum_{k=0}^{\infty}\mathbb{P}[X
\geq k].
\]
With $X=\frac1 {\varepsilon}\log(1+|\xi_0|)$ it follows that
$\mathbb{E} \log(1+|\xi_0|)<\infty$ if and only if $
\sum_{k=1}^{\infty}\mathbb{P}[|\xi_0|\geq e^{\varepsilon k}-1] <
\infty%
$ for some (equivalently, every) $\varepsilon>0$. The proof is
completed by
applying the Borel--Cantelli lemma.
\end{pf}

Note in passing that Lemma~\ref{lemlogmoment} and condition~(A4) imply
that for every $n\in\mathbb{N}$ the series~(\ref{eqdefGGG})
converges with
probability $1$ on $\mathbb{D}_{R_0}$.

\subsection{Upper bound in Theorem~\texorpdfstring{\protect\ref{theologasympt}}{4.1}}
\label{subsecupperbound}
Fix an $\varepsilon>0$. All constants which we will introduce below
depend only on $\varepsilon$. Let us agree that all inequalities will
hold uniformly over $z\in\mathbb{D}_{e^{-2\varepsilon}
R_0}\setminus\{0\}$ if $R_0<\infty$ and over
$z\in\mathbb{D}_{1/\varepsilon}\setminus\{0\}$ if $R_0=\infty $. We
will show that there exists an a.s. finite random variable
$M=M(\varepsilon)$ such that
for all sufficiently large $n$, 
%
\begin{equation}
\label{equpperboundG} \bigl|\mathbf{G}_n(z) \bigr|\leq M e^{n(I(\log|z|)+3\varepsilon)}.
\end{equation}
First, we estimate the tail of the Taylor series~(\ref{eqdefGGG})
defining $\mathbf{G}_n$. 
By assumption~(A4) there is $A>\max(0,-\log f(0))$ such that for all
$n\geq A$ and all $k\geq An$,
\[
|f_{k,n}|< \bigl(|z| e^{2\varepsilon} \bigr)^{-k}.
\]
Lemma~\ref{lemlogmoment} implies that there exist a.s. finite random
variables $S,M'$ such that for all $n\geq A$,
%
\begin{equation}\label{equpperbound1}
\biggl|\sum_{k\geq An}\xi_k f_{k,n} z^k \biggr|\leq S \sum
_{k\geq An} e^{\varepsilon k} |f_{k,n}|
|z|^k \leq S \sum_{k\geq An} e^{-\varepsilon k}
\leq M' e^{-A n}.
\end{equation}
%

We now consider the initial part of the Taylor series~(\ref{eqdefGGG})
defining $\mathbf{G}_n$. Take some $\delta>0$. By assumption~(A3),
there is $N$ such that for all $n>N$ and all $k\leq An$,
%
\begin{equation}
\label{eqestfkn} |f_{k,n}|< \biggl(f \biggl(\frac kn \biggr)+\delta
\biggr)^n.
\end{equation}
It follows from~(\ref{eqdefI}) that for all $t\geq0$,
%
\begin{equation}
\label{eqestIlogz} t\log|z| + \log f(t) \leq I\bigl(\log|z|\bigr).
\end{equation}
Using~(\ref{eqestfkn}), (\ref{eqestIlogz}) and Lemma~\ref{lemlogmoment}
with $\varepsilon/A$ instead of $\varepsilon$ we obtain that there is
an a.s. finite
random variable $M''$ such that for all sufficiently large $n$,
%
\begin{eqnarray}
\label{equpperbound2}\quad
\biggl\vert\sum_{0\leq k< An}\xi_k f_{k,n} z^k \biggr\vert
&\leq& M'' \sum_{0\leq k< An}
e^{(\varepsilon k)/A} \biggl(f \biggl(\frac kn \biggr)+ \delta \biggr)^n
|z|^k
\\
&\leq& M'' e^{\varepsilon n} \sum_{0\leq k< An} \bigl(e^{(k/n)\log|z|+\log
f(k/n)}+ \delta|z|^{k/n}\bigr)^n
\nonumber
\\
&\leq& M'' e^{2 \varepsilon n} \bigl(e^{I(\log|z|)}+
\delta\max \bigl(1, |z|^A \bigr) \bigr)^n
\nonumber
\\
&\leq& M'' e^{n (I(\log|z|)+3\varepsilon)},
\nonumber
\end{eqnarray}
where the last inequality holds if $\delta=\delta(\varepsilon)$ is
sufficiently small. Combining~(\ref{equpperbound1})
and~(\ref{equpperbound2}) and noting that $-A<\log f(0)\leq I(\log
|z|)$ by~(\ref{eqestIlogz}), we obtain that~(\ref{equpperboundG}) holds
with $M=M'+M''$ for sufficiently large $n$. By enlarging $M$, if
necessary, we can achieve that it holds for all $n\geq A$.


\subsection{Lower bound in Theorem~\texorpdfstring{\protect\ref{theologasympt}}{4.1}}\label{subseclowerbound}
Fix $\varepsilon>0$ and $z\in\mathbb{D}_{R_0}\setminus\{0\}$. We are
going to show that
%
\begin{equation}
\label{eqprooflogpartfunc2} \mathbb{P} \bigl[\bigl|\mathbf{G}_n(z)\bigr|<
e^{n(I(\log|z|)-4\varepsilon)} \bigr]=O \biggl(\frac1 {\sqrt n} \biggr),\qquad n\to\infty.
\end{equation}
We will use the Kolmogorov--Rogozin inequality in a multidimensional
form which can be found in~\cite{esseen}. Given a $d$-dimensional
random vector $X$ define its concentration function by
%
\begin{equation}
\label{eqdefconcentrationfunc} Q(X; r) = \sup_{x\in\mathbb{R}^d}
\mathbb{P} \bigl[X\in\mathbb{D}_r(x) \bigr], \qquad r>0,
\end{equation}
where $\mathbb{D}_r(x)$ is a $d$-dimensional ball of radius $r$
centered at $x$. An easy consequence of~(\ref{eqdefconcentrationfunc})
is that for all independent random vectors $X,Y$ and all $r,a>0$,
%
\begin{equation}
\label{eqpropconcentrationfunc} Q(X+Y; r)\leq Q(X; r), \qquad Q(aX;r)=Q(X; r/a).
\end{equation}
The next result follows from Corollary~1 on page~304 of~\cite{esseen}.

%
\begin{theorem}[(Kolmogorov--Rogozin inequality)]
\label{theokolmogorov-rogozin} There is a constant $C_d$ depending only
on $d$ such that for all independent (not necessarily identically
distributed) random $d$-dimensional vectors $X_1,\ldots,X_n$ and for
all $r>0$, we have
\[
Q(X_1+\cdots+X_n; r) \leq C_d \cdot
\Biggl( \sum_{k=1}^n \bigl(1-Q(X_k;
r) \bigr) \Biggr)^{-1/2}.
\]
\end{theorem}

The idea of our proof of~(\ref{eqprooflogpartfunc2}) is to use the
Kolmogorov--Rogozin inequality to show that the probability of very
strong cancellation among the terms of the series~(\ref{eqdefGGG})
defining $\mathbf{G}_n$ is small. First, we have to single out those
terms of $\mathbf{G}_n$ in which $|f_{k,n}z^k|$ is large enough. By
definition of $I$, see~(\ref{eqdefI}), there is $t_0\in[0, T_0]$ such
that $t_0\log |z|+\log f(t_0)>I(\log|z|)-\varepsilon$. Moreover, by
assumption~(A2), we can find a closed interval $J$ of length $|J|>0$
containing $t_0$ such that
\[
f(t) |z|^t >e^{I(\log|z|)-2\varepsilon},\qquad t\in J.
\]
%
Define a set $\mathcal{J}_n=\{k\in\mathbb{N}_0\dvtx  k/n \in J\}$. By
assumption~(A3) there is $N$ such that for all $n>N$ and all
$k\in\mathcal{J}_n$,
\[
|f_{k,n}| |z|^k > e^{n(I(\log|z|)-3\varepsilon)}.
\]
Let $n>N$. For $k\in\mathbb{N}_0$ define
\[
a_{k,n}=e^{-n(I(\log|z|)-3\varepsilon)} f_{k,n} z^k.
\]
Note that $|a_{k,n}| > 1$ for $k\in\mathcal{J}_n$. Define
\[
\mathbf{G}_{n,1} = \sum_{k\in\mathcal{J}_n}
a_{k,n} \xi_k,\qquad\mathbf{G}_{n,2} = \sum
_{k\notin\mathcal{J}_n} a_{k,n}\xi_k.
\]
By considering real and imaginary parts, we can view the complex random
variables $a_{k,n}\xi_k$ as two-dimensional random vectors.
Using~(\ref{eqpropconcentrationfunc}), we arrive at
%
\begin{equation}\label{eqprooflogpartfunc2a}
\qquad\mathbb{P} \bigl[\bigl|\mathbf{G}_n(z)\bigr|<
e^{n(I(\log|z|)-4\varepsilon
)} \bigr] \leq Q \bigl( \mathbf{G}_{n,1}+
\mathbf{G}_{n,2}; e^{-\varepsilon n} \bigr) \leq Q \bigl(
\mathbf{G}_{n,1}; e^{-\varepsilon n} \bigr).
\end{equation}
By Theorem~\ref{theokolmogorov-rogozin}, there is an absolute constant
$C$ such that for all $r>0$,
\begin{eqnarray*}
Q(\mathbf{G}_{n,1}; r) &\leq& C\cdot \biggl(\sum
_{k\in\mathcal{J}_n} \bigl(1-Q(a_{k,n}\xi_k; r) \bigr)
\biggr)^{-1/2}
\\
&\leq& C\cdot \biggl(\sum_{k\in\mathcal{J}_n}
\bigl(1-Q( \xi_k; r) \bigr) \biggr)^{-1/2}.
\end{eqnarray*}
Here, the second inequality follows from the fact that $|a_{k,n}| > 1$
for $k\in\mathcal{J}_n$. Now, since the random variable $\xi_0$ is
supposed to be nondegenerate, we can choose $r>0$ so small that
$Q(\xi_0; r)<1$. Note that this is the only place in the proof of
Theorem~\ref{theogeneral} where we use the randomness of the $\xi_k$'s
in a nonobvious way. The rest of the proof is valid for any
deterministic sequence $\xi_0,\xi_1,\ldots$ such that
$|\xi_n|=O(e^{\delta n})$ for every $\delta>0$. If $n$ is sufficiently
large, then $e^{-\varepsilon n}\leq r$ and hence,
%
\begin{equation}
\label{eqprooflogpartfunc2b} Q \bigl(\mathbf{G}_{n,1};
e^{-\varepsilon n} \bigr)\leq Q(\mathbf{G}_{n,1}; r) \leq
C_1 |\mathcal{J}_n|^{-1/2} \leq C_2
n^{-1/2}.
\end{equation}
In the last inequality, we have used that the number of elements of
$\mathcal{J}_n$ is larger than $(|J|/2) n$ for large $n$.
Taking~(\ref{eqprooflogpartfunc2a}) and~(\ref{eqprooflogpartfunc2b})
together completes the proof of~(\ref{eqprooflogpartfunc2}).

\subsection{Proof of Proposition~\texorpdfstring{\protect\ref{lemconvpotconvmeas}}{4.3}}
\label{subsecconvpotconvmeas} Define a set
$A\subset\mathbb{D}_{R_0}\times\Omega$, measurable with respect to the
product of the Borel $\sigma$-algebra on $\mathbb{D}_{R_0}$ and
$\mathcal{F}$, by
\[
A= \Bigl\{(z, \omega)\dvtx \lim_{k\to\infty} p_{l_k}(z;
\omega) = I\bigl(\log|z|\bigr) \Bigr\}. 
\]
We know from assumption~(\ref{eqplktoI}) that for Lebesgue-a.e.
$z\in\mathbb{D}_{R_0}$ it holds that $\int_{\Omega} \mathbh{1}_{(z,
\omega)\notin
A}\mathbb{P}(d\omega)=0$. By Fubini's theorem, for $\mathbb{P}$-a.e.
$\omega\in
\Omega$,
it holds that $\int_{\mathbb{D}_{R_0}} \mathbh{1}_{(z, \omega
)\notin A}
\lambda(d
z)=0$. Hence, there is a measurable set $E_1\subset\Omega$ with
$\mathbb{P}[E_1]=0$ such that for every $\omega\notin E_1$,
%
\begin{equation}
\label{eqlimitplk} \lim_{k\to\infty} p_{l_k}(z;\omega) =
I\bigl(\log|z|\bigr),\qquad\mbox{for Lebesgue-a.e. }z\in\mathbb{D}_{R_0}.
\end{equation}

Let $k(\omega)=\min\{k\in\mathbb{N}_0\dvtx  \xi_k(\omega)\neq0\}$,
$\omega\in\Omega$. Since the $\xi_k$'s are assumed to be nondegenerate,
the set $E_0=\{\omega\in\Omega\dvtx  k(\omega)=\infty\}$ satisfies
$\mathbb{P}[E_0]=0$. By conditions~(A3) and~(A1), after ignoring
finitely many values of $n$, we can assume that $f_{k,n}\neq0$ for
$0\leq k \leq T_0n/2$. Define $n(\omega)=2k(\omega)/T_0$. For $\omega
\notin E_0$ and $n>n(\omega)$ the function $\mathbf{G}_n$ does not
vanish identically. For every fixed $\omega\notin E_0$ and
$n>n(\omega)$ the function $p_n(z;\omega)=\frac1n
\log|\mathbf{G}_n(z;\omega)|$ is subharmonic, as a function of $z$; see
Example 4.1.10 in~\cite{hoermanderbook1}. Also, it follows
from~(\ref{equpperboundG}) that there is a measurable set
$E_2\subset\Omega$ with $\mathbb{P}[E_2]=0$ such that for every
$\omega\notin E_2$, the family of functions $\mathcal{P}_{\omega}=\{ z
\mapsto p_{l_k}(z;\omega)\dvtx  k\in\mathbb{N}\}$, is uniformly bounded
above on every compact subset of $\mathbb{D}_{R_0}$. Let $E=E_0\cup
E_1\cup E_2$, so that $\mathbb{P}[E]=0$. Fix $\omega\notin E$. By
Theorem~4.1.9 of~\cite{hoermanderbook1}, the family
$\mathcal{P}_{\omega}$ is either precompact in
$\mathcal{D}'(\mathbb{D}_{R_0})$, the space of generalized functions on
the disk $\mathbb{D}_{R_0}$, or contains a subsequence converging to
$-\infty$ uniformly on compact subsets of $\mathbb{D}_{R_0}$. The
latter possibility is excluded by~(\ref{eqlimitplk}). Thus, the family
$\mathcal{P}_{\omega}$ is precompact in $\mathcal{D}'(\mathbb
{D}_{R_0})$. Any subsequential limit of $\mathcal{P}_{\omega}$ must
coincide with the function $I(\log|z|)$ by~(\ref{eqlimitplk}) and
Proposition~16.1.2 in~\cite{hoermanderbook2}. It follows that for every
fixed $\omega\notin E$,
%
\begin{equation}
\label{eqlimitplk1} p_{l_k}(z; \omega) \mathop{\longrightarrow}_{k\to\infty}I\bigl(
\log |z|\bigr)\qquad\mbox{in }\mathcal{D}'( \mathbb{D}_{R_0}).
\end{equation}
%
Since the Laplace operator is continuous on $\mathcal{D}'(\mathbb
{D}_{R_0})$, we
may apply it to the both sides of~(\ref{eqlimitplk1}).
Recalling~(\ref{eqmuGGGisLaplacelogGGG}) we obtain that for every
$\omega\notin E$,
\[
\frac{1}{l_k}\mu_{l_k}(\omega)= \frac{1}{2\pi} \Delta
p_{l_k}(z; \omega) \mathop{\longrightarrow}_{k\to\infty}
\frac
{1}{2\pi} \Delta I\bigl(\log|z|\bigr)\qquad \mbox{in }\mathcal{D}'(
\mathbb{D}_{R_0}).
\]
A sequence of locally finite measures converges in $\mathcal
{D}'(\mathbb{D}_{R_0})$
if and only if it converges vaguely. This completes the proof
of~(\ref{eqlimitmulk}).

\subsection{Proof of the a.s. convergence in Theorem~\texorpdfstring{\protect\ref{theoLOpoly}}{2.3}} \label{subsecproofLOpolyas}
Recall that convergence in
probability has already been established in
Section~\ref{secproofsspecialcases}. To prove the a.s. convergence we
first extract a subsequence to which we can apply the Borel--Cantelli
lemma. Given $n\in\mathbb{N}$ we can find a unique $j_n\in\mathbb {N}$
such that $j_n^3\leq n< (j_n+1)^3$. Write $m_n=j_n^3$ and
$\mathbf{G}_n(z)=\mathbf{W}_n(e^{\beta} m_n^{\alpha}z)$. Note that
$\lim_{n\to\infty}m_n/n =1$. Thus, it suffices to show that $\frac1n
\mu_{\mathbf{G}_n}$ converges a.s. to the measure with
density~(\ref{eqLOpolydensity}). As a first step, we will prove the
a.s. convergence of the corresponding potentials. Fix $z\in\mathbb
{D}\setminus \{0\}$. We will prove that
%
\begin{equation}
\label{eqLOasconvpotentials} 
p_n(z) = \frac1n \log\bigl|
\mathbf{G}_n(z)\bigr| \mathop{\longrightarrow }^{\mathrm{a.s.}}_{n\to\infty}
\alpha|z|^{1/\alpha}.
\end{equation}
Note that $\mathbf{G}_n$ satisfies all assumptions of
Section~\ref{subsecgeneraltheo}. It follows from
Proposition~\ref{proplogasymptsubseq} applied to the subsequence
$l_j=j^3$ that
%
\begin{equation}
\label{eqLOasconvpotentialssubseq} \frac1{m_n} \log\bigl|
\mathbf{G}_{m_n}(z)\bigr| \mathop{\longrightarrow }^{\mathrm{a.s.}}_{n\to\infty}
\alpha|z|^{1/\alpha}.
\end{equation}
%
Let now $n\in\mathbb{N}$ be a sufficiently large number not of the
form $j^3$.
We have, by Lemma~\ref{lemlogmoment} and~(\ref{eqfkasympt}),
\begin{eqnarray*}
\bigl|\mathbf{G}_n(z) - \mathbf{G}_{m_n}(z)\bigr| &=& \Biggl\vert\sum
_{k=m_n+1}^{n} \xi_k w_k e^{\beta k} m_n^{\alpha k} z^k\Biggr\vert 
\\
&\leq& S e^{2\varepsilon n} \sum_{k=m_n+1}^{n} e^{-\alpha(k\log k-k)} n^{\alpha k}
|z|^k.
\end{eqnarray*}
The function $x\mapsto-\alpha(x\log x-x)+\alpha x \log n$ defined for
$x>0$ attains its maximum, which is equal to $\alpha n$, at $x=n$.
Recall that $|z|<1$. Since $m_n>(1-\varepsilon)n$ and $\varepsilon n
S<e^{\varepsilon n}$ if
$n$ is sufficiently large, we have the estimate
\[
\bigl|\mathbf{G}_n(z)-\mathbf{G}_{m_n}(z)\bigr| \leq
e^{3\varepsilon n} e^{\alpha n} |z|^{(1-\varepsilon) n}.
\]
Since $\alpha+\log|z| < \alpha|z|^{1/\alpha}$, we have, if
$\varepsilon>0$
is small enough,
%
\begin{equation}
\label{eqasconvlowerdiff} \bigl| \mathbf{G}_n(z)-\mathbf{G}_{m_n}(z)\bigr|
\leq e^{(1-\varepsilon)n(\alpha
|z|^{1/\alpha
}-2\varepsilon
)}\leq e^{m_n (\alpha|z|^{1/\alpha}-2\varepsilon)}.
\end{equation}
Bringing~(\ref{eqLOasconvpotentialssubseq})
and~(\ref{eqasconvlowerdiff}) together we
obtain~(\ref{eqLOasconvpotentials}).

We are ready to complete the proof. It follows
from~(\ref{eqLOasconvpotentials}) and
Proposition~\ref{lemconvpotconvmeas} that the restriction of $\frac1n
\mu_{\mathbf{G}_n}$ to $\mathbb{D}$ converges a.s. to a measure $\mu $
with density~(\ref{eqLOpolydensity}), as a sequence of random elements
with values in $\mathcal{M}(\mathbb{D})$. To prove that the a.s.
convergence holds in the sense of $\mathcal{M}(\mathbb{C})$-valued
elements, we need to show that $\lim_{n\to\infty} \frac1n
\mu_{\mathbf{G}_n}(\mathbb {C}\setminus \mathbb{D})=0$ a.s., or,
equivalently, that $\liminf_{n\to\infty} \frac1n \mu_{\mathbf{G}
_n}(\mathbb{D})=1$ a.s. Let $f\dvtx \mathbb{C}\to[0,1]$ be a continuous
function with support in $\mathbb{D}$. Then, since \mbox{$\nu\mapsto\int f
\,d\nu$} defines a continuous functional on $\mathcal{M}(\mathbb{D})$,
\[
\liminf_{n\to\infty} \frac1n \mu_{\mathbf{G}_n}(\mathbb{D}) \geq
\liminf_{n\to
\infty} \frac1n \int_{\mathbb{C}} f\,d
\mu_{\mathbf{G}_n} = \int_{\mathbb{C}} f\,d\mu\qquad\mbox{a.s.}
\]
The supremum of the right-hand side over all admissible $f$ is equal to
$1$ since $\mu(\mathbb{D})=1$. This proves the claim.
\subsection{Proof of the a.s. convergence in Theorem~\texorpdfstring{\protect\ref{theoLOentire}}{2.5}}
\label{subsecproofLOentireas} Let $m_n$ be defined in
the same way as in the previous proof. Write $\mathbf{G}_n(z)=\mathbf
{W}(e^{\beta}
m_n^{\alpha}z)$. Note that $\mathbf{G}_n$ satisfies the assumptions of
Section~\ref{subsecgeneraltheo} with $I(s)=\alpha e^{s/\alpha}$.
By Proposition~\ref{proplogasymptsubseq}, for all $z\in\mathbb
{C}\setminus\{0\}$,
%
\begin{equation}
\label{eqLOasconvpotentialsentire} 
p_n(z) = \frac1n
\log\bigl| \mathbf{G}_n(z)\bigr| \mathop{\longrightarrow }^{\mathrm{a.s.}}_{n\to\infty}
\alpha|z|^{1/\alpha}.
\end{equation}
Then, it follows from Proposition~\ref{lemconvpotconvmeas} that
$\frac1n \mu_{\mathbf{G}_n}$ converges a.s. to the measure with
density~(\ref{eqLOentiredensity}).


\subsection{Proof of Theorem~\texorpdfstring{\protect\ref{theoLOpolyconverse}}{2.4}}
We prove only the implication $(1)\Rightarrow(2)$ since the converse
implication has been established in Theorem~\ref{theoLOpoly}. Let
$\mathbf{W}_n(z)=\sum_{k=0}^n \xi_k w_k z^k$, where $w_k$ is a sequence
satisfying~(\ref{eqfkasympt}) and~(\ref{eqLOadditionalassumpt}).
Assume that $\mathbb{E}\log(1+|\xi_0|)=\infty$. Fix $\varepsilon
>0$. We will show that
with probability $1$ there exist infinitely many $n$'s such that all
zeros of $\mathbf{W}_n(e^{\beta}n^{\alpha} z)$ are located in the disk~$\mathbb{D}_{2\varepsilon}$. This implies that $\frac1n \mu_n$ does
not converge
a.s. to the measure with density~(\ref{eqLOpolydensity}). We use an
idea of~\cite{izlog}. By Lemma~\ref{lemlogmoment},
$\limsup_{n\to\infty} |\xi_n|^{1/n}=+\infty$. Hence, with probability~$1$ there exist infinitely many $n$'s such that
%
\begin{eqnarray}\label{eqnologmomentsuperexp}
|\xi_n|^{1/n}&>&\max_{k=1,\ldots,n-1} | \xi_{n-k}|^{1/(n-k)},
\nonumber\\[-8pt]\\[-8pt]
|\xi_n|^{1/n}&>&\max \biggl\{\frac{3C+1}{\varepsilon}, \frac1
{e^{\alpha}\varepsilon} \biggr\}.\nonumber
\end{eqnarray}
Let $n$ be such that~(\ref{eqnologmomentsuperexp}) holds.
By~(\ref{eqLOadditionalassumpt}) and~(\ref{eqnologmomentsuperexp}), we
have for every $z\in\mathbb{C}$ and $k<n$,
\begin{eqnarray*}
\bigl\vert w_{n-k}\xi_{n-k} \bigl(e^{\beta} n^{\alpha}
z \bigr)^{n-k} \bigr\vert&\leq& C |w_n| e^{\beta k}n^{\alpha k}
|\xi_n|^{(n-k)/n} \bigl| e^{\beta} n^{\alpha} z
\bigr|^{n-k}
\\
&=& C \bigl\vert w_n \xi_n \bigl(e^{\beta}n^{\alpha}z
\bigr)^n \bigr\vert \bigl(|\xi_n|^{1/n} |z|
\bigr)^{-k}.
\end{eqnarray*}
For every $z$ such that $|z|>\varepsilon$, we obtain
\begin{eqnarray*}
\Biggl\vert\sum_{k=1}^{n-1} w_{n-k}
\xi_{n-k} \bigl(e^{\beta} n^{\alpha
}z \bigr)^{n-k}\Biggr\vert&\leq& C \bigl\vert w_n \xi_n \bigl(e^{\beta}
n^{\alpha}z \bigr)^n \bigr\vert\cdot \Biggl(\sum
_{k=1}^{n-1} \frac1{(3C+1)^{k}} \Biggr)
\\
&<& \frac13 \bigl\vert w_n \xi_n \bigl(e^{\beta}
n^{\alpha}z \bigr)^n \bigr\vert.
\end{eqnarray*}
By~(\ref{eqfkasympt}) and~(\ref{eqnologmomentsuperexp}), the right-hand
side of this inequality goes to $+\infty$ as $n\to\infty$. In
particular, for sufficiently large $n$, it is larger than $|\xi_0w_0|$.
It follows that for $|z|>\varepsilon$, the term of degree $n$ in the
polynomial $\mathbf{W}_n(e^{\beta}n^{\alpha} z)$ is larger, in the
sense of absolute value, than the sum of all other terms. Hence, the
polynomial $\mathbf{W}_n(e^{\beta}n^{\alpha} z)$ has no zeros outside
the disk $\mathbb{D}_{2\varepsilon}$.

\subsection{Proof of Theorem~\texorpdfstring{\protect\ref{theogeneralconverse}}{2.9}} 
Start with a measure $\mu$ satisfying the assumptions of
Theorem~\ref{theogeneralconverse}. Define a function $I$ by
$I(s)=\int_{-\infty}^{s}\mu(\mathbb{D}_{e^r})\,dr$ for $s<\log R_0$.
The integral is finite by the second assumption of the theorem.
Clearly, $I$ is nondecreasing, continuous and convex on $(-\infty, \log
R_0)$. For $s>\log R_0$ let $I(s)=+\infty$. Define $I(\log R_0)$ by
left continuity. Let now $u$ be defined as the Legendre--Fenchel
transform of $I$:
\[
u(t)=\sup_{s\in\mathbb{R}} \bigl(st- I(s) \bigr).
\]
We claim that the random analytic function
$\mathbf{G}_n(z)=\sum_{k=0}^{\infty} \xi_k f_{k,n} z^k$ with
$f_{k,n}=e^{-n u(k/n)}$ satisfies assumptions (A1)--(A4) of
Theorem~\ref{theogeneral} with $f=e^{-u}$. By the Legendre--Fenchel
duality, the function $u$ possesses the following properties. First, it
is convex and lower-semicontinuous. Second, it is finite on the
interval $[0,T_0)$, where $T_0=\limsup_{t\to+\infty} I(t)/t$ satisfies
$T_0\in (0,+\infty]$. This holds since $I$ is nondecreasing and
$\lim_{s\to -\infty} I(s)=0$ by construction. Third, $u(t)=+\infty$ for
$t>T_0$ and $t<0$. This verifies assumption~(A1). Fourth,
formula~(\ref{eqdefI}) holds and $\lim_{t\to+\infty} u(t)/t=\log R_0$.
This, together with Lemma~\ref{lemlogmoment}, shows that the
convergence radius of $\mathbf{G}_n$ is $R_0$ a.s. and verifies
assumption~(A4). Finally, $u$ is continuous on $[0,T_0)$ (since it is
convex and finite there), and, in the case $T_0<+\infty$,\vadjust{\goodbreak} the function
$u$ is left continuous at $T_0$ (follows from the lower-semicontinuity
of $u$). This verifies assumption~(A2). Assumption~(A3) holds trivially
with $f=e^{-u}$.

\section*{Acknowledgment}
The authors are grateful to the unknown referee who considerably
simplified the original proof of Theorem~\ref{theogeneral}. The
argument in Section~\ref{subsecconvpotconvmeas} follows the idea of the
referee.



%

\printaddresses

\end{document}